\documentclass[12pt]{article}
\usepackage{amsmath,amsfonts,latexsym,amsthm,amssymb}
\newcommand{\K}{\overline{K}}
\DeclareMathOperator{\card}{card}

\DeclareMathOperator{\Dim}{Dim}
\DeclareMathOperator{\const}{const}
\begin{document}
\newtheorem{prop}{Proposition}
\newtheorem{teo}{Theorem}
\pagestyle{plain}
\title{Analysis and Probability over Infinite Extensions of a Local Field,
II: A Multiplicative Theory}

\bigskip
\author{ANATOLY N. KOCHUBEI\footnote{E-mail: ank@ank.kiev.ua}
\quad Institute of Mathematics,
National Academy\\ of Sciences of Ukraine, Tereshchenkivska 3,
Kiev, 01601 Ukraine}
\date{}
\maketitle
\begin{abstract}
Let $V$ be a projective limit, with respect to the renormalized
norm mappings, of the groups of principal units corresponding to
a strictly increasing sequence of finite separable totally and
tamely ramified Galois extensions of a local field. We study the
structure of the dual group $V'$, introduce and investigate a
fractional differentiation operator on $V$, and the corresponding
L\'evy process.
\end{abstract}
\section{INTRODUCTION}

This paper is a continuation of the article \cite{K1} (see also
\cite{K2,Y}), in which we considered an infinite extension $K$ of
a local field of zero characteristic assuming that $K$ is a union
of an increasing sequence of finite extensions. If $K$ is
equipped with a natural inductive limit topology, its strong
conjugate $\K$ is a projective limit with respect to the
renormalized trace mappings. A Gaussian measure, a Fourier
transform, a fractional differentiation operator, and a cadlag
Markov process $X_\alpha$ (an analog of the $\alpha$-stable
process) on $\K$ were constructed. The semigroup of measures
defining $X_\alpha$ is concentrated on a compact (additive)
subgroup $S\subset \K$. Sample paths properties of the part of
$X_\alpha$ in $S$ were studied in \cite{K3}.

All the above constructions were based essentially on algebraic
structures related to the additive groups of the field $K$ and
its subfields. In this paper we develop a parallel theory based
on the multiplicative structures. Note that in analysis over
local fields both the approaches are closely connected (see
Sections 3.5 and 4.7 in \cite{K2}).

Let $k$ be a non-Archimedean local field of an arbitrary
characteristic. We consider a strictly increasing sequence of its
finite separable Galois extensions
\begin{equation}
k=K_1\subset K_2\subset \ldots \subset K_n\subset \ldots .
\end{equation}
We shall assume that all the extensions (1) are totally and
tamely ramified. Denote $m_n=[K_n:K_1]$, $n=2,3,\ldots$. It will
be convenient to write $m_0=0$, $m_1=1$. By our assumptions, the
residue field cardinality for each field $K_n$ is the same
positive integer $q=p^\varkappa$ where $p$ is the characteristic
of the residue fields. If $\nu >n$, then $[K_\nu
:K_n]=\dfrac{m_\nu}{m_n}$.

Let $U_{1,n}$ be the group of principal units of the field $K_n$.
Then
\begin{equation*}
U_{1,1}\subset U_{1,2}\subset \ldots \subset U_{1,n}\subset \ldots .
\end{equation*}
Under the above assumptions the norm mapping $N_{\nu ,n}:\ K_\nu \to
K_n$ maps $U_{1,\nu}$ onto $U_{1,n}$ ($\nu >n$); see Chapter 1,
\S 8 in \cite{ANT}. For any $x\in U_{1,n}$ denote
$$
\mathfrak N_{\nu ,n}(x)=\left[ N_{\nu ,
n}(x)\right]^{\frac{m_n}{m_\nu}}.
$$
Since $\dfrac{m_\nu}{m_n}$ is prime to $p$, $\mathfrak N_{\nu ,
n}(x)$ is well-defined and belongs to $U_{1,n}$ (\cite{FV},
Chapter 1, Corollary (5.5)). It is easy to check that $\mathfrak
N_{\nu ,n}(x):\ U_{1,\nu}\to U_{1,n}$ is an epimorphism, and
$\mathfrak N_{\nu ,n}(x)=x$ if $x\in U_{1,n}$.

If $l>\nu >n$, then for any $x\in U_{1,l}$
$$
\mathfrak N_{\nu ,n}\left( \mathfrak N_{l,\nu }(x)\right) =\left[
N_{\nu ,n}\left( \left( N_{l,\nu
}(x)\right)^{\frac{m_\nu}{m_l}}\right)\right]^{\frac{m_n}{m_\nu}}.
$$
Note that for any positive integer $r$ \ $N_{\nu ,n}(x^r)=\left[
N_{\nu ,n}(x)\right]^r$. If $r$ is prime to $p$, we may substitute
$x^{1/r}$ for $x$ and find that $N_{\nu ,n}(x)=\left[
N_{\nu ,n}(x^{1/r})\right]^r$, so that $N_{\nu ,n}$ commutes with
the $r$-th root operation. Therefore $\mathfrak N_{\nu ,n}\circ
\mathfrak N_{l,\nu }=\mathfrak N_{l,n}$, and we can define the
projective limit
$$
V=\varprojlim U_{1,n}
$$
with respect to the homomorphisms $\mathfrak N_{\nu ,n}$.

$V$ is a compact topological group. Its structure is investigated
in Sect. 2. $V$ is totally disconnected, its topology is
determined by an explicitly written descending chain of open-closed
subgroups. The structure of its dual group is studied.

In Sect. 3 we introduce a fractional differentiation operator
$D^\alpha$ on $V$ (continuous characters on $V$ are its
eigenfunctions), and prove that it is a generator of a Markov
process $\xi_\alpha$ on $V$, an analog of the stable process.
Sect. 4 is devoted to properties of $\xi_\alpha$. In particular,
we find the Hausdorff and packing dimensions of the image of a
time interval under $\xi_\alpha$. This is based on the explicit
calculation of the L\'evy measure corresponding to $\xi_\alpha$,
and employs the general results by Evans \cite{E} on sample paths
properties of L\'evy processes on Vilenkin groups.

\section{THE GROUP $V$ AND ITS DUAL}

2.1. The group $V$ consists of sequences $x=(x_1,\ldots
,x_n,\ldots )$, $x_n\in U_{1,n}$, such that for any $\nu >n$
\ $\mathfrak N_{\nu ,n}(x_\nu )=x_n$, with the component-wise
operations. The topology on $V$ is induced by the Tychonoff
topology on the direct product $\prod\limits_nU_{1,n}$. Thus a
fundamental system of neighbourhoods of the unit element in $V$
can be obtained by taking finite intersections of the sets
$$
\left\{ x=(x_1,\ldots ,x_n,\ldots )\in V:\ |1-x_n|_n\le
\varepsilon \right\}
$$
with some $n\ge 1$, $0<\varepsilon <1$. Here $|\cdot |_n$ is the
normalized absolute value on $K_n$. Note that each of the above
sets is an open-closed subgroup in $V$ due to the ultra-metric
property of the absolute values.

Below we shall use the following property of the mappings
$\mathfrak N_{\nu ,n}$, $\nu >n$. Let
$$
U_{l,i}=\left\{ y\in K_i:\ |1-y|_i\le q^{-l}\right\},\quad l\ge
2.
$$
Then
\begin{equation}
\mathfrak N_{\nu ,n}:\ U_{\frac{m_\nu}{m_n}(l-1)+1,\nu
}\xrightarrow{\text{onto}}U_{l,n}\ ,
\end{equation}
\begin{equation}
\mathfrak N_{\nu ,n}:\ U_{\frac{m_\nu}{m_n}(l-1),\nu
}\xrightarrow{\text{onto}}U_{l-1,n}\ .
\end{equation}
For the mappings $N_{\nu ,n}$ this is proved in \cite{S}, Chapter
V, \S 6, Corollary 3, where a more general case is considered;
the Hasse-Herbrand function $\psi (v)$ of an extension $L/K$,
which appears in that corollary, equals $[L:K]v$ for a totally and
tamely ramified extension (see \cite{S}, Chapter IV, \S 3). In order
to consider the mappings $\mathfrak N_{\nu ,n}$ it remains to use the
$m$-divisibility of the groups $U_{l,n}$ for any $m$ prime to $p$
(\cite{FV}, Chapter 1, (5.5)).

Let us consider subgroups $V_n\subset V$,
$$
V_n=\left\{ x=(x_1,\ldots ,x_n,\ldots )\in V:\ |1-x_n|_n\le
q^{-nm_n-1} \right\},\quad n=1,2,\ldots .
$$
If $x\in V_n$ and $n>i$, then $x_i=\mathfrak N_{n,i}(x_n)$, and by
(2)
$$
|1-x_i|_i\le q^{-m_i(n+1)}<q^{-im_i-1},
$$
so that $x_i\in V_i$. Thus we have a filtration
\begin{equation}
V=V_0\supset V_1\supset \ldots \supset V_n\supset \ldots .
\end{equation}
The same arguments show that the subgroups $V_n$ form a
fundamental system of neighbourhoods of the unit element.

Note that the descending chain (4) can be ``lengthened'' by including
intermediate subgroups so that the resulting chain would be such
that the quotient group of two consequtive subgroups is of a
prime order. This property is often assumed in the investigation
of totally disconnected groups. In particular, this assumption
was made in \cite{E}. However it will be more convenient for us
to use the chain (4). All the results of \cite{E} remain valid
here.

Let $\mu$ be the normalized Haar measure on $V$. Since the
mappings $\mathfrak N_{\nu ,n}$ are surjective, $V$ can be
identified with the projective limit of its quotient groups with
respect to a descending system of compact subgroups (\cite{B1},
Chapter III, \S 7). Then $\mu$ is a projective limit of the
normalized Haar measures $\mu_n$ on the groups $U_{1,n}$,
$n=1,2,\ldots$ (\cite{B2}, Chapter VII, \S 1, Sect. 6). This
means in particular that
$$
\mu \left( V_n\right) =\mu_n \left( U_{nm_n+1,n}\right) .
$$

On the other hand, $\mu_n$ is proportional to the normalized Haar
measure $dx_n$ on the additive group of the field $K_n$. Since
(see e.g. \cite{K2}) $\int\limits_{U_{1,n}}dx_n=q^{-1}$, we have
$$
\mu_n \left( U_{nm_n+1,n}\right) =q\int\limits_{|1-x_n|_n\le q^{-
nm_n-1}}=q^{-nm_n}.
$$

It is known \cite{V,E} that $V/V_n$ is a finite group of the
order $M(n)$ where $\mu \left( V_n\right) =[M(n)]^{-1}$. From the
above calculations we find that
\begin{equation}
M(n)=q^{nm_n}.
\end{equation}

\medskip
2.2. By the duality theorem \cite{L}, the dual discrete group
$V'$ is the inductive limit $\varinjlim U_{1,n}$ with respect to
the dual mappings $\mathfrak N'_{\nu ,n}:\ U'_{1,n}\to U'_{1,\nu}$
($\nu >n$) defined as follows. If $\theta_n\in U'_{1,n}$, that is
$\theta_n$ is a (multiplicative) continuous character of the group
$U_{1,n}$, then  $\mathfrak N'_{\nu ,n}(\theta_n)$ is a character
of the group $U_{1,\nu}$, and for any $x_\nu \in U_{1,\nu}$
$$
\mathfrak N'_{\nu n}(\theta_n)(x_\nu )=\theta_n\left( \mathfrak N_{\nu
n}(x_\nu)\right).
$$

By definition (see \cite{L} or \cite{F}), $V'=\mathfrak
A/\mathfrak B$ where $\mathfrak A$ is the direct sum
$\bigoplus\limits_{n=1}^\infty U'_{1,n}$, $\mathfrak B$ is a
subgroup consisting of all elements $\left( \theta_{i_1},\ldots
,\theta_{i_n}\right)$, $\theta_{i_j}\in U'_{1,j}$, such that for
some $i\ge i_1,\ldots ,i_n$
$$
\left( \mathfrak N'_{i,i_1}\theta_{i_1}\right) \cdots
\left( \mathfrak N'_{i,i_n}\theta_{i_n}\right) =1;
$$
in our multiplicative notation we identify $\left( \theta_{i_1},\ldots
,\theta_{i_n}\right)$ with $\theta_{i_1}\ldots \theta_{i_n}$.
Each element of $V'$ can be written as $\theta =\theta_i\mathfrak
B$, $\theta_i\in U_{1,i}'$ for some $i$. Then the coupling
between $V'$ and $V$ is given by $\langle \theta ,x\rangle =
\langle \theta_i ,x_i\rangle_i$, if $x=(x_1,\ldots ,x_i,\ldots
)$. Here $\langle \cdot ,\cdot \rangle_i$ is the coupling between
$U'_{1,i}$ and $U_{1,i}$. The correctness of this definition is
proved in \cite{L}.

In our specific situation we can say a little more about the
above definition of $V'$. Suppose that $\theta_i^{(1)}\in
U'_{1,i}$, $\theta_j^{(2)}\in U'_{1,j}$, and
\begin{equation}
\mathfrak N'_{n,i}\theta_i^{(1)}=\mathfrak
N'_{n,j}\theta_j^{(2)},\quad i\le j\le n.
\end{equation}
The equality (6) means that
$$
\theta_i^{(1)}\left( \mathfrak N_{n,i}(x)\right) =\theta_j^{(2)}\left(
\mathfrak N_{n,j}(x)\right)\quad \mbox{for any \ }x\in U_{1,n}.
$$
In particular, taking $x\in U_{1,j}\subset U_{1,n}$ we find that
$\mathfrak N_{n,j}(x)=x$, $\mathfrak N_{n,i}(x)=\mathfrak
N_{j,i}\left( \mathfrak N_{n,j}(x)\right) =\mathfrak N_{j,i}(x)$,
so that
$$
\mathfrak N'_{j,i}\theta_i^{(1)}=\theta_j^{(2)},\ \mbox{if
}i<j;\quad \theta_i^{(1)}=\theta_j^{(2)},\ \mbox{if }i=j.
$$

Thus, two representatives of the same coset either coincide (if
they belong to the same group $U'_{1,j}$), or one of them is
obtained by ``lifting'' another.

In order to classify elements of $V'$, we shall need the notion
of a ramification degree of a multiplicative character of a local
field (adapted to our setting). Let $\theta_n\in U'_{1,n}$. The
character $\theta_n$ is said to have the ramification degree 1,
if $\theta_n\equiv 1$, and the ramification degree $\nu \ge 2$,
if $\theta_n(x_n)=1$ for any $x_n\in U_{\nu,n}$, and
$\theta_n(x^0_n)\ne 1$ for some $x^0_n\in U_{\nu -1,n}$.

It follows from (2) and (3) that if $\theta_i\in U'_{1,i}$ has
the ramification degree $\nu_i$, and $n>i$, then $\mathfrak
N'_{n,i}\theta_i$ has the ramification degree
$\nu_n=\dfrac{m_n}{m_i}(\nu_i-1)+1$. We see that
$$
\frac{\nu_n-1}{m_n}=\frac{\nu_i-1}{m_i}.
$$
Therefore we may assign to any $\theta \in V'$ the number
$$
r(\theta )=\frac{\nu_i-1}{m_i}
$$
where $\nu_i$ is the ramification degree of an arbitrary
representative of $\theta$ lying in $U'_{1,i}$.

\medskip
2.3. Let $V_n^\bot\subset V'$ be the annihilator of the subgroup
$V_n$. We have
$$
\{1\}=V_0^\bot\subset V_1^\bot \subset \ldots \subset V_n^\bot
\subset \ldots ;\quad \bigcup\limits_{n=0}^\infty V_n^\bot =V'.
$$
It is known \cite{V,E} that $\card \left( V_n^\bot\right) =M(n)$.

\begin{prop}
The annihilator $V_n^\bot$ consists of those cosets $\theta \in
V'$ for which $r(\theta )\le n$, and there exists a representative
$\theta_n\in \theta$, $\theta_n\in U'_{1,n}$.
\end{prop}

{\it Proof}. Let us show first of all that a coset $\theta \in
V_n^\bot$ contains a representative $\theta_n\in U'_{1,n}$.
Indeed, if $\theta_i\in \theta$, and $i<n$, then the character
$\theta_i$ can be lifted (within the coset $\theta$) to the
character $\theta_n(x_n)=\theta_i\left( \mathfrak
N_{n,i}(x_n)\right)$, $x_n\in U_{1,n}$.

Suppose that $i>n$. Let us define a character $\theta_n\in V_n'$
as a restriction of $\theta_i$ to $U_{1,n}$. We have to check
that $\theta_i\theta_n^{-1}\in \mathfrak B$; it is sufficient to
verify that $\theta_i\cdot \left( \mathfrak
N'_{i,n}\theta_n\right) =1$, that is
$$
\theta_i(x_i)\left[ \theta_n\left( \mathfrak N_{i,n}(x_i)\right)
\right]^{-1}=1
$$
for any $x_i\in U_{1,i}$, or (by the definition of $\theta_n$)
that for each $x_i\in U_{1,i}$
\begin{equation}
\theta_i\left( \frac{x_i}{\mathfrak N_{i,n}(x_i)}\right) =1.
\end{equation}

Denote $\widetilde{x}_i=\dfrac{x_i}{\mathfrak N_{i,n}(x_i)}$,
$\widetilde{x}_l=\mathfrak N_{i,l}\left( \widetilde{x}_i\right)$
for $l<i$ (since $\mathfrak N_{i,n}(x_i)\in U_{1,n}\subset
U_{1,i}$, the element $\widetilde{x}_i$ belongs to $U_{1,i}$).
Choosing $\widetilde{x}_{i+1}\in U_{1,i+1}$ in such a way that
$\mathfrak N_{i+1,i}\left( \widetilde{x}_{i+1}\right) =\widetilde{x}_i$,
then taking such $\widetilde{x}_{i+2}\in U_{1,i+2}$ that
$\mathfrak N_{i+2,i+1}\left( \widetilde{x}_{i+2}\right) =\widetilde{x}_{i+1}$
etc., we obtain an element $\widetilde{x}=\left( \widetilde{x}_1,\ldots
,\widetilde{x}_n,\ldots ,\widetilde{x}_i,\ldots \right) \in V$. We
have
$$
\widetilde{x}_n=\left\{ \frac{N_{i,n}(x_i)}{\left[ \mathfrak
N_{i,n}(x_i)\right]^{\frac{m_i}{m_n}}}\right\}^{\frac{m_n}{m_i}}=1,
$$
so that $\widetilde{x}\in V_n$, which means that $1=\theta
(\widetilde{x})=\theta_i\left( \widetilde{x}_i\right)$, and we have
proved the equality (7).

Thus, $\theta =\theta_n\cdot \mathfrak B$, $\theta_n\in
U'_{1,n}$. If $r(\theta )>n$, then the ramification degree of the
character $\theta_n$ is greater than $nm_n+1$, which implies the
existence of an element $x\in V_n$, such that $\theta (x)\ne 1$
(as in the above reasoning, an element $x_n\in U_{nm_n+1,n}$ can
be prolonged to an element $(x_1,\ldots ,x_n,\ldots )\in V_n$).
This means that the inequality $r(\theta )>n$ implies that
$\theta \notin V_n^\bot$. On the other hand, if $r(\theta )\le n$
and $\theta =\theta_n\cdot \mathfrak B$, $\theta_n\in U'_{1,n}$,
then obviously $\theta \in V_n^\bot$. $\quad \blacksquare$

\section{Fractional Differentiation Operator}

3.1. Let us consider the function $f$ on $V'$ defined as
$$
f(\theta )=\varphi_n\ \ \mbox{for }\theta \in V_n^\bot\setminus
V_{n-1}^\bot ,\ n=1,2,\ldots ;\ f(1)=\varphi_0,
$$
where $\{ \varphi_n\}_{n=0}^\infty$ is a given sequence of
complex numbers. It follows from (5) that $f\in l_1(V')$, if
$\sum\limits_{n=0}^\infty |\varphi_n|q^{nm_n}<\infty$, and $f\in l_2(V')$, if
$\sum\limits_{n=0}^\infty |\varphi_n|^2q^{nm_n}<\infty$.

Let us compute the Fourier transform
$$
F(x)=\sum\limits_{\theta \in V'}f(\theta )\theta (x),\quad x\in
V.
$$
We have
$$
F(x)=\varphi_0+\sum\limits_{n=1}^\infty \varphi_n
\sum\limits_{\theta \in V_n^\bot\setminus V_{n-
1}^\bot }\theta (x).
$$
If $x\in V\setminus V_1$, then $\sum\limits_{\theta \in
V_n^\bot}\theta (x)=0$ for $n\ge 1$, whence $F(x)=\varphi_0-
\varphi_1$. If $x\in V_l\setminus V_{l+1}$, $l\ge 1$, then
$$
F(x)=\varphi_0+\sum\limits_{n=1}^l\varphi_n\sum\limits_{\theta
\in V_n^\bot\setminus V_{n-1}^\bot }1
+\varphi_{l+1}\left[ \sum\limits_{\theta \in V_{l+1}^\bot }\theta
(x)-\sum\limits_{\theta \in V_l^\bot }1\right] +\sum\limits_{n=l+2}^\infty
\varphi_n\sum\limits_{\theta \in V_n^\bot\setminus V_{n-
1}^\bot }\theta (x).
$$
Since $\sum\limits_{\theta \in V_{l+1}^\bot }\theta (x)=0$ and
$\sum\limits_{\theta \in V_n^\bot\setminus V_{n-1}^\bot }\theta
(x)=0$ for $n\ge l+2$, and $\sum\limits_{\theta \in V_n^\bot
}1=q^{nm_n}$, we find that
$$
F(x)=\varphi_0+\sum\limits_{n=1}^l\varphi_n\left[ q^{nm_n}-q^{(n-
1)m_{n-1}}\right] -\varphi_{l+1}q^{lm_l}.
$$
This gives after a simple transformation that
\begin{equation}
F(x)=\sum\limits_{n=0}^l(\varphi_n-\varphi_{n+1})q^{nm_n},\quad
x\in V_l\setminus V_{l+1},\ l\ge 0.
\end{equation}

If $\sum\limits_{n=0}^\infty |\varphi_n|q^{nm_n}<\infty$, then it
follows from (8) that
$$
F(1)=\sum\limits_{n=0}^\infty (\varphi_n-\varphi_{n+1})q^{nm_n}.
$$

\medskip
3.2. Consider the function $f^{(\alpha )}(\theta )$, $\theta \in
V'$, corresponding, as above, to the sequence
$$
\varphi_n^{(\alpha )}=\begin{cases}
q^{\alpha nm_n}, & \text{if $n\ge 1$},\\
0, & \text{if $n=0$},\end{cases}
$$
where $\alpha <-1$. Let $F^{(\alpha )}$ be the Fourier transform
of $f^{(\alpha )}$. Introducing the convolution operator
$A^{(\alpha )}u=F^{(\alpha )}*u$, $u\in L_2(V)$, we can write
$$
\left( A^{(\alpha )}u\right) (x)=\int\limits_VF^{(\alpha
)}(y)u(xy^{-1})\mu (dy)=\int\limits_VF^{(\alpha
)}(y)\left[ u(xy^{-1})-u(x)\right] \mu (dy),
$$
since
$$
\int\limits_VF^{(\alpha )}(y)\mu (dy)=f^{(\alpha
)}(1)=\varphi_0^{(\alpha )}=0.
$$

Denote by $\mathcal D(V)$ the vector space of locally constant
complex-valued functions on $V$, that is such functions $u$ that
$u(x)=u(y)$ if $xy^{-1}\in V_l$ (the number $l$ depends on $u$
and does not depend on $x$). Since $\bigcup_{n=0}^\infty V_n^\bot
=V'$, any continuous character on $V$ belongs to $\mathcal D(V)$.
Therefore $\mathcal D(V)$ is dense in the Banach space $C(V)$ of
all continuous functions on $V$.

By (8),
\begin{equation}
F^{(\alpha )}(y)=-q^\alpha +\sum\limits_{n=1}^l\left( q^{\alpha
nm_n}-q^{\alpha (n+1)m_{n+1}}\right) q^{nm_n},\quad y\in
V_l\setminus V_{l+1},
\end{equation}
(the sum is missing for $l=0$). If $u\in \mathcal D(V)$, then
$\left( A^{(\alpha )}u\right) (x)$ is an entire function of
$\alpha$, and we define the operator $D^\alpha u$, $\alpha >0$,
as the analytic continuation of $A^{(\alpha )}u$. Thus,
$$
\left( D^\alpha u\right) (x)=\int\limits_VF^{(\alpha
)}(y)\left[ u(xy^{-1})-u(x)\right] \mu (dy),\quad \alpha >0,
$$
for any $u\in \mathcal D(V)$. The expression (9) is valid for
$\alpha >0$ too. Equivalently, $D^\alpha$ can be written on
$\mathcal D(V)$ as a pseudo-differential operator with the symbol
$f^{(\alpha )}(\theta )$.

\medskip
\begin{teo}
{\rm (i)} The operator $D^\alpha$ ($\alpha >0$) is an essentially
selfadjoint operator on $L_2(V)$. Its closure has a purely discrete
spectrum consisting of the eigenvalues
$$
\varphi_n^{(\alpha )}=\begin{cases}
q^{\alpha nm_n}, & \text{if $n\ge 1$},\\
0, & \text{if $n=0$},\end{cases}
$$
corresponding to the eigenspaces $V_n^\bot\setminus V_{n-1}^\bot
$, $n=0,1,2,\ldots$ ($V_{-1}^\bot =\varnothing$). In particular,
$D^\alpha 1=0$.

{\rm (ii)}  The semigroup $e^{-tD^\alpha }$, $t\ge 0$, consists
of the integral operators of the form
$$
\left( e^{-tD^\alpha }u\right) (x)=\int\limits_VG_\alpha (t,xy^{-
1})u(y)\mu (dy)
$$
where the kernel
\begin{equation}
G_\alpha (t,z)=\sum\limits_{n=0}^l\left[ e^{-t\varphi_n^{(\alpha
)}}-e^{-t\varphi_{n+1}^{(\alpha )}}\right] q^{nm_n},\quad z\in
V_l\setminus V_{l+1},\ l\ge 0,
\end{equation}
is positive. The corresponding stochastic process $\xi_\alpha$ with
independent increments on $V$ is stochastically continuous.
\end{teo}

{\it Proof}. All the assertions except the last one are immediate
consequences  of the above constructions. Note that the function
$z\mapsto G_\alpha (t,z)$ is constant on each open-closed set
$V_l\setminus V_{l+1}$, which implies its continuity. Thus
$\xi_\alpha$ has a strong Feller property. The stochastic
continuity is equivalent to the $C_0$-property of the semigroup
$e^{-tD^\alpha }$ in $C(V)$. It is sufficient to prove that
$\left\| e^{-tD^\alpha }u-u\right\|_{C(V)}\to 0$ as $t\to 0$, for
any $u\in \mathcal D(V)$.

Since $\int\limits_VG_\alpha (t,xy^{-1})\mu (dy)=1$, we find that
for any $u\in \mathcal D(V)$
\begin{multline*}
\left( e^{-tD^\alpha }u-u\right) (x)=\int\limits_VG_\alpha (t,xy^{-
1})[u(y)-u(x)]\mu (dy)=\int\limits_VG_\alpha (t,z)u[(xz^{-1})-u(x)]\mu (dz)
\\ =\int\limits_{V\setminus V_l}G_\alpha (t,z)u[(xz^{-1})-u(x)]\mu (dz)
=\sum\limits_{j=0}^{l-1}\int\limits_{V_j\setminus V_{j+1}}G_\alpha
(t,z)u[(xz^{-1})-u(x)]\mu (dz)
\end{multline*}
for some $l$. Now it follows from (10) that
$$
\left\| e^{-tD^\alpha }u-u\right\|_{C(V)}\le \const \cdot
\|u\|_{C(V)}\sum\limits_{j=0}^{l-1}\sum\limits_{n=0}^j
\left[ e^{-t\varphi_n^{(\alpha
)}}-e^{-t\varphi_{n+1}^{(\alpha )}}\right]
q^{nm_n}\longrightarrow 0,\quad t\to 0,
$$
as desired. $\quad \blacksquare$

\section{Sample Path Properties}

4.1. It is well known (see e.g. \cite{H}) that for any $\theta
\in V'$
\begin{equation}
\mathbf E\langle \theta ,\xi_\alpha (t)\rangle =\exp \left\{
\int\limits_V[\langle \theta ,x\rangle -1]\Pi (t,dx)\right\}
\end{equation}
where $\Pi$ is the L\'evy measure, that is $\Pi (t,\Gamma
)=\mathbf EN(t,\Gamma )$, $N(t,\Gamma )$ is the number of jumps
of $\xi_\alpha$ on the interval $[0,t)$ belonging to a Borel set
$\Gamma \not\ni 1$.

\begin{prop}
The L\'evy measure, $\Pi (t,\Gamma )$ has the form $\Pi (t,\Gamma
)=t\nu (\Gamma )$ where $\nu (\Gamma )$ is a Borel measure on
$V\setminus \{1\}$ finite outside any open neighbourhood of 1,
which possesses the following properties:

{\rm (i)} $\nu$ satisfies the condition of local spherical
symmetry \cite{E}, that is $\nu (dx)=\nu_n\mu (dx)$ on
$V_n\setminus V_{n+1}$. Specifically,
\begin{equation}
\nu_n=q^\alpha +\sum\limits_{l=1}^nq^{lm_l}\left[ q^{\alpha
(l+1)m_{l+1}}-q^{\alpha lm_l}\right] .
\end{equation}

{\rm (ii)} The equality
\begin{equation}
\nu (V\setminus V_n)=\sum\limits_{j=0}^{n-1}\nu_j\left[ q^{-
jm_j}-q^{-(j+1)m_{j+1}}\right]
\end{equation}
holds.

{\rm (iii)} The asymptotic relation
\begin{equation}
\nu (V\setminus V_n)\sim q^{\alpha nm_n},\quad \mbox{as }n\to
\infty ,
\end{equation}
is valid.
\end{prop}

\medskip
{\it Proof}. By the construction of the process $\xi_\alpha$,
$$
\mathbf E\langle \theta ,\xi_\alpha (t)\rangle
=\int\limits_V\langle \theta ,z\rangle G_\alpha (t,z)\mu
(dz)=e^{-tf^{(\alpha )}(\theta )}
=\begin{cases}
e^{-tq^{\alpha nm_n}}, & \text{if $\theta \in V_n^\bot \setminus
V_{n-1}^\bot ,n\ge 1$},\\
1, & \text{if $\theta \equiv 1$}.\end{cases}
$$
Since $G_\alpha (t,z)=G_\alpha (t,z^{-1})$, the measure $\Pi$ is
also invariant with respect to the inversion $z\mapsto z^{-1}$.
Therefore comparing the last inequality with (11) we find that
$$
\int\limits_V[\langle \theta ,x\rangle -1]\Pi (t,dx)=\begin{cases}
-tq^{\alpha nm_n}, & \text{if $\theta \in V_n^\bot \setminus
V_{n-1}^\bot ,n\ge 1$},\\
0, & \text{if $\theta \equiv 1$}.\end{cases}
$$

Let $\Gamma \subset V_n\setminus V_{n+1}$ be a Borel set,
$\omega_\Gamma (x)$ the indicator of the set $\Gamma$. We have
$$
\omega_\Gamma (x)=\sum\limits_{\theta \in V'}\widehat{\omega_\Gamma}
(\theta )\theta (x),\quad \widehat{\omega_\Gamma}
(\theta )=\int\limits_\Gamma \overline{\langle \theta
,x\rangle}\mu (dx).
$$
In particular, $\sum\limits_{\theta \in V'}\widehat{\omega_\Gamma}
(\theta )=\omega_\Gamma (1)=0$, so that
$$
\omega_\Gamma (x)=\sum\limits_{\theta \in V'}\widehat{\omega_\Gamma}
(\theta )[\theta (x)-1],
$$
whence
\begin{multline*}
\int\limits_V\omega_\Gamma (x)\Pi (t,dx)=-
t\sum\limits_{l=1}^\infty q^{\alpha lm_l}\sum\limits_{\theta \in
V_l^\bot \setminus V_{l-1}^\bot}\int\limits_\Gamma \overline{\langle \theta
,x\rangle}\mu (dx)=-t\sum\limits_{l=1}^n q^{\alpha
lm_l}\card \left( V_l^\bot \setminus V_{l-1}^\bot \right) \mu (\Gamma )\\
-tq^{\alpha (n+1)m_{n+1}}\left[ \int\limits_\Gamma \left(
\sum\limits_{\theta \in V_{n+1}^\bot}\langle \theta ,x\rangle \right) \mu
(dx)-\mu (\Gamma )\card V_n^\bot \right]\\ -
t\sum\limits_{l=n+2}^\infty q^{\alpha lm_l}\left\{ \int\limits_\Gamma
\left[ \sum\limits_{\theta \in V_l^\bot }\langle \theta ,x\rangle
\right] \mu (dx)-\int\limits_\Gamma
\left[ \sum\limits_{\theta \in V_{l-1}^\bot }\langle \theta ,x\rangle
\right] \mu (dx)\right\} .
\end{multline*}

An element $x\in V_n\setminus V_{n+1}$ can be considered as a
non-trivial character of each group $V_l^\bot$, $l\ge n+1$ (since
$\left( V_l^\bot \right)^\bot =V_l$). Therefore $\sum\limits_{\theta
\in V_l^\bot }\langle \theta ,x\rangle =0$ for $l\ge n+1$, so
that
\begin{multline*}
\int\limits_V\omega_\Gamma (x)\Pi (t,dx)=-t\mu (\Gamma )\left\{
\sum\limits_{l=1}^n q^{\alpha lm_l}\left[ q^{lm_l}-q^{(l-1)m_{l-
1}}\right] -q^{\alpha (n+1)m_{n+1}}\cdot q^{nm_n}\right\} \\
=t\mu (\Gamma)\left\{
\sum\limits_{l=1}^n q^{lm_l}\left[ q^{\alpha (l+1)m_{l+1}}-
q^{\alpha lm_l}\right] +q^\alpha \right\} ,
\end{multline*}
and we have come to (12), which easily implies (13).

In order to prove (14), we apply the Abel transform to the sum in
(12). We have
\begin{equation}
\nu_n=q^\alpha +q^{nm_n}\left( q^{\alpha (n+1)m_{n+1}}-q^\alpha
\right) -\sum\limits_{i=1}^{n-1}\left( q^{(i+1)m_{i+1}}-
q^{im_i}\right) \left( q^{\alpha (i+1)m_{i+1}}-q^\alpha
\right) .
\end{equation}
Since $m_{n+1}\ge 2m_n$, it follows from (15) that
\begin{equation}
\nu_n\sim q^{{nm_n}+\alpha (n+1)m_{n+1}},\quad n\to \infty .
\end{equation}

Now the asymptotics (14) is a consequence of (13) and (16).
Indeed, the right-hand side of (14) is clearly the ``leading''
term in (13). As for other summands of (13), we find that
$$
\nu_{n-1}q^{-nm_n}=\left( \nu_{n-1}q^{-(n-1)m_{n-1}}\right) \cdot
q^{(n-1)m_{n-1}-nm_n}\le \const \cdot q^{\alpha nm_n}\cdot q^{-
\left( \frac{n}2+1\right) m_n}=o\left( q^{\alpha nm_n}\right) ;
$$
$$
\sum\limits_{j=0}^{n-2}\nu_j\left[ q^{-jm_j}-q^{-
(j+1)m_{j+1}}\right] \le \const \cdot \sum\limits_{j=0}^{n-
2}q^{\alpha (j+1)m_{j+1}}\le \const \cdot (n-
1)q^{\frac{\alpha}2(n-1)m_n}=o\left( q^{\alpha nm_n}\right)
$$
because $m_n\ge 2^{n-1}$. $\quad \blacksquare$

\medskip
4.2. Using Proposition 2 and the general results by Evans
\cite{E}, we can find now the Hausdorff and packing dimensions of
the image of a time interval under the random mapping $t\mapsto
\xi_\alpha (t)$. See \cite{E} for the definitions. Below we
assume that $\xi_\alpha (0)=1$.

Let $\pi (n)$ be the first exit time of the process $\xi_\alpha$
out of the subgroup $V_n$. Let
$$
Q(n,N)=\mathbf P\left\{ \xi_\alpha (t)\notin V_n \ \forall t\in
[\pi (n),\pi (N))\right\} ,\quad n>N.
$$
An essential assumption in \cite{E} is the inequality
\begin{equation}
\liminf\limits_{n\to \infty }Q(n,N)>0.
\end{equation}
In specific situations the verification of (17) can be a
difficult task; see, for example, \cite{K3}. However if the
L\'evy measure is locally spherically symmetric, as we have here
by Proposition 2, then (17) takes place if
\begin{equation}
\limsup\limits_{n\to \infty }\frac{\nu (V\setminus
V_n)}{M(n)}\cdot \frac{M(n-1)}{\nu (V\setminus V_{n-1})}<1
\end{equation}
(\cite{E}, Corollary 2).

It follows from (5) and (14) that the inequality (18) holds if
$\alpha <1$.

Denote by $\dim$ the Hausdorff dimension, and by $\Dim$ the
packing dimension corresponding to the natural metric on $V$ (see
\cite{E}).

\begin{teo}
If $\alpha <1$, then for each $t>0$ we have that
\begin{equation}
\dim \xi_\alpha ([0,t])=\Dim \xi_\alpha ([0,t])=\alpha
\end{equation}
almost surely.
\end{teo}

{\it Proof}. The equalities (19) follow from (17) and Theorems 8,
10 of \cite{E}. $\quad \blacksquare$

\medskip
Note that the Hausdorff and packing dimensions of an image of an
interval for a stable process on the field $\mathbb Q_p$ of $p$-adic
numbers were found recently in \cite{AZ}.

\end{document}